\documentclass[12pt]{article}
\usepackage{graphicx,subeqn}
\usepackage{pstricks}
\usepackage{amsfonts,amssymb}
\usepackage[numbers]{natbib}
\bibpunct{(}{)}{;}{a}{,}{,}
\usepackage[bookmarksnumbered=true, pdfauthor={Wen-Long Jin}]{hyperref}

\newcommand{\commentout}[1]{}

\newcommand{\ba}{\begin{array}}
        \newcommand{\ea}{\end{array}}
\newcommand{\bc}{\begin{center}}
        \newcommand{\ec}{\end{center}}
\newcommand{\bdm}{\begin{displaymath}}
        \newcommand{\edm}{\end{displaymath}}
\newcommand{\bds} {\begin{description}}
        \newcommand{\eds} {\end{description}}
\newcommand{\ben}{\begin{enumerate}}
        \newcommand{\een}{\end{enumerate}}
\newcommand{\beq}{\begin{equation}}
        \newcommand{\eeq}{\end{equation}}
\newcommand{\bfg} {\begin{figure}[p]}
        \newcommand{\efg} {\end{figure}}
\newcommand{\bi} {\begin {itemize}}
        \newcommand{\ei} {\end {itemize}}
\newcommand{\bpp}{\begin{pspicture}}
        \newcommand{\epp}{\end{pspicture}}
\newcommand{\bqn}{\begin{eqnarray}}
        \newcommand{\eqn}{\end{eqnarray}}
\newcommand{\bqs}{\begin{eqnarray*}}
        \newcommand{\eqs}{\end{eqnarray*}}
\newcommand{\bsq}{\begin{subequations}}
        \newcommand{\esq}{\end{subequations}}
\newcommand{\bsl} {\begin{slide}[8.8in,6.7in]}
        \newcommand{\esl} {\end{slide}}
\newcommand{\bss} {\begin{slide*}[9.3in,6.7in]}
        \newcommand{\ess} {\end{slide*}}
\newcommand{\btb} {\begin {table}[p]}
        \newcommand{\etb} {\end {table}}
\newcommand{\bvb}{\begin{verbatim}}
        \newcommand{\evb}{\end{verbatim}}

\newcommand{\m}{\mbox}
\newcommand {\der}[2] {{\frac {\m {d} {#1}} {\m{d} {#2}}}}
\newcommand {\pd}[2] {{\frac {\partial {#1}} {\partial {#2}}}}

\newcommand{\cas}[1]{{{\left \{ \ba #1 \ea \right. }}}

\newcommand{\reff}[1] {{{Figure \ref {#1}}}}
\newcommand{\refe}[1] {{{(\ref {#1})}}}
\newcommand{\reft}[1] {{{\textbf{Table} \ref {#1}}}}

\newtheorem{theorem}{Theorem}[section]
\newtheorem{definition}[theorem]{Definition}
\newtheorem{lemma}{Lemma}

\def\pmb#1{\setbox0=\hbox{$#1$}%
   \kern-.025em\copy0\kern-\wd0
   \kern.05em\copy0\kern-\wd0
   \kern-.025em\raise.0433em\box0 }

\def\eop{{\hfill $\blacksquare$}}

\def\J {{\bf J}}
\def\eop{{\hfill $\blacksquare$}}

\title{A dynamical system model of the traffic assignment problem}
\author{Wen-Long Jin \thanks{Institute of Transportation Studies,  University of California, 522 Social Science Tower, Irvine, CA 92697, USA. Email: wjin@uci.edu. Corresponding author}}
\begin {document}
\maketitle
\begin{abstract}
User equilibrium is a central concept for studying transportation networks, and one can view it as the result of a dynamical process of drivers' route choice behavior.
In this paper, based on a definition of O-D First-In-First-Out violation, we propose a new dynamical system model of the route choice behavior at the aggregate, route level for both static and dynamic transportation networks. An equilibrium of such a dynamical system can be a user equilibrium or a partial user equilibrium. We prove that, for static, symmetric traffic assignment problem with fixed or variable demand, only user equilibria are stable for the dynamical system, and the objective function in the  mathematical programming formulation \citep*{Beckmann1956bmw} can be considered as the potential energy of the dynamical system. We then present an Euler-based perturbation method for finding user equilibrium and solve two examples for both static and dynamic traffic assignment problems. This new model is simple in form and could be applied to analyze other properties of transportation networks.
\end{abstract}


\section{Introduction}
The traffic assignment problem (TAP) is to assign {\em origin-destination} (O-D) flows to a road network. Solving this problem is a core step of four-step transportation planning method and essential for analyzing congestion and other performance of a network. It is generally believed that assigned traffic flows are in the state of {\em user equilibrium} (UE) \citep{wardrop1952ue}, where ``the journey times of all routes actually used are equal and less than those which would be experienced by a single vehicle on any unused route". That is, for the same O-D pair, all used routes share the same travel time, and unused ones are not shorter. A fundamental task in solving TAP is to compute travel times on all links. In static TAP, it is assumed that there exist so-called {\em link performance functions}, which map link flows to travel times. For static TAP, link performance functions can be separable, symmetric, or monotone. While in dynamic TAP, not only the number of vehicles but the dynamics of traffic can affect the travel time on a link, and there may not exist any meaningful link performance functions \citep{daganzo1995tt}.

There are four well-known formulations of TAP: namely, mathematical programming, variational inequality, nonlinear complementarity, and fixed point formulations  \citep{Patriksson1994ta}. In addition, there exists another so-called {\em dynamical system} approach, including dynamical models of stochastic, individual route choice behavior \citep[e.g.][]{Horowitz1984stability} and dynamical models of aggregate, deterministic change of route flows \citep{smith1984stability,Nagurney1996ds,Nagurney1997equilibria}. Equilibria of these dynamical systems are usually UE. However, these dynamical systems are very complex in form and have not been widely used in solving and analyzing static or dynamic UE for large-scale transportation networks.

In this paper, we propose a new dynamical system formulation of TAP. For static, symmetric link performance functions, we define a new aggregate, deterministic route choice dynamics based on a so-called {\em First-In-First-Out (FIFO) violation} function among different routes of the same O-D pair. We then show that, in equilibria of this dynamical system, ``the journey times of all routes actually used are equal", but may not be ``less than those which would be experienced by a single vehicle on any unused route". That is, an equilibrium of this dynamical system may be a UE or a {\em partial user equilibrium} (PUE). For static, symmetric TAP, however, we are able to show that PUE are always unstable, and only UE are stable. Then, based on the stability properties of UE and PUE, we propose a numerical algorithm for computing stable UE. We further extend this new dynamical system formulation for TAP with variable demand and for dynamic TAP that incorporates traffic dynamics. Through this study, we intend to develop a model that can be used for theoretical analysis of properties of various transportation networks.

The rest of the paper is organized as follows. In Section 2, we introduce a new dynamical system of static TAP and analyze stability properties of its equilibria. In Section 3, we present a finite difference method and a perturbation-based method for computing equilibria of the dynamical system and UE. In Section 4, we extend the new dynamical system formulation for dynamic TAP. We make some concluding remarks in Section 5.

\section{A dynamical system model of static traffic assignment problem}
For static TAP, we adopt the network notation system given in \reft{notation.2}, which is similar to that in \citep{sheffi1985networks}, but with several new definitions. We also have the following basic relationships:
\bsq\label{conservation}
\bqn
q_{rs}=\sum_k f_k^{rs}, && \forall\:\: r,s \\
f_k^{rs}\geq 0, && \forall\:\: k,r,s; \label{nonnegative}
\eqn
\esq
and
\bsq\label{basic}
\bqn
x_a=\sum_{r,s}\sum_k f_k^{rs} \delta_{a,k}^{rs}, && \forall\:\: a\\
c_k^{rs}=\sum_a t_a \delta_{a,k}^{rs}, && \forall\:\: k,r,s.
\eqn
\esq

We call ${\bf f}$ as a UE state in the sense of \citep{wardrop1952ue} if and only if for ${\bf f}\in \mathcal{F}$
\bsq
\bqn
f_k^{rs}(c_k^{rs}-v_{rs})=0, && \forall\:\: k,r,s\\
c_k^{rs}\geq v_{rs}, && \forall\:\: k,r,s.
\eqn
\esq
As we know, UE is the solution of the following mathematic programming problem \citep*{Beckmann1956bmw}
\bqn
\min z(\bf x)&=&\sum_a \int_0^{x_a} t_a(\omega) d\omega, \label{objective}
\eqn
whose objective function is hereafter called the BMW objective function. The BMW objective function $z(\bf x)$ is strictly convex in $\bf x$ for static, symmetric TAP with fixed demand \citep{Dafermos1969ta}. In this section, we assume that $z(\bf x)$ is convex, but not necessarily strictly convex. That is, $z({\bf x})$ still attains its minimum at UE solutions of link flows, but UE solutions may not be unique in ${\bf x}$. The BMW objective function can also be written in route flows, $z({\bf f})\equiv z({\bf x}({\bf f}))$, and, at UE, $z({\bf f})$ also attains its minimum.

\subsection{O-D FIFO violation and a dynamical system model of aggregate route choice behavior}
On a road link, FIFO among vehicles is violated when one vehicle passes another. Such link FIFO violation can be caused by heterogeneity between vehicles, different traffic conditions on different lanes, or traffic signals. Measurements were developed for link FIFO violation among different groups of vehicles in \citep{jin2005fifo} and among individual vehicles in \citep{jin2006fifo}. For vehicles of the same O-D, even when vehicles observe the FIFO principle on all links and routes, the FIFO principle can still be violated if vehicles using different routes experience different travel times. Here we are interested in such {\em O-D FIFO violation} and assume that in O-D FIFO state vehicles departing from the origin at the same time arrive at the destination at the same time. Hereafter, FIFO violation means O-D FIFO violation if not otherwise noted.

To show how to measure O-D FIFO violation in static transportation networks, we consider a network with one O-D pair, $r-s$, and two alternative routes. In one assignment, $(f_1^{rs}, f_2^{rs})$, one obtains two different route travel times, $(c_1^{rs}, c_2^{rs})$. The cumulative arrival curves at the origin and destination \citep{moskowitz1963capacity,newell1993sim} are shown in \reff{fifostatic}, where $v_{rs}$ is the average O-D travel time. Then we can define FIFO violation for two route flows, $(J_1^{rs}, J_2^{rs})$, by the areas of the two regions shown in \reff{fifostatic}, but subject to opposite signs. Note that here the unit of FIFO violation is vehicle$\times$time.

For a general road network, we can define the FIFO violation function for route $k$ connecting O-D pair $r-s$ by ($\forall\:\: k, r,s$)
\bsq
\bqn
J_k^{rs}&=&q_{rs}f_k^{rs}(c_k^{rs}-v_{rs}), \label{odfv}
\eqn
whose unit is vehicle$^2\times$time \footnote{Here we multiply FIFO violation defined in \reff{fifostatic} by $q_{rs}$ for convenience, and our discussions apply if all O-D FIFO violation functions of an O-D pair are multiplied by the same constant.}. Since the average O-D travel time is
\bqn
v_{rs}&=&\frac{\sum_j f_j^{rs} c_j^{rs}}{q_{rs}},\label{averageOD}
\eqn
\esq
the FIFO violation function can be re-written as
\bqn
J_k^{rs}&=&f_k^{rs}(q_{rs}c_k^{rs}-\sum_j f_j^{rs} c_j^{rs})=f_k^{rs}\sum_j f_j^{rs} (c_k^{rs}- c_j^{rs}).
\eqn

To model the route choice dynamics at the aggregate, route level, we introduce a new principle so that the traffic flow on a route is decreased when the travel time on the route is longer than the average O-D travel time. In particular, assuming that all route flows depend on a decision variable $\tau$, whose unit is $(vehicle\times time)^{-1}$, we propose that the rate of change in a route flow be negative to the corresponding O-D FIFO violation; that is,
\bqn
-\dot f_k^{rs}=J_k^{rs}=f_k^{rs}\sum_j (c_k^{rs}-c_j^{rs})f_j^{rs}=q_{rs}f_k^{rs}(c_k^{rs}-v_{rs}), && \forall\:\: k,r,s \label{jsystem:1}
\eqn
where $\dot f_k^{rs}$ is the derivative of $f_k^{rs}$ with respect to $\tau$.
From \refe{jsystem:1}, we can have the following autonomous dynamical system of route flows ${\bf f}$
\bqn
-\dot {{\bf f}}&=&\J({\bf f}).  \label{fifodyn}
\eqn
Note that, with the definition of O-D FIFO violation in \refe{odfv}, we have a negative sign on the left-hand side of \refe{fifodyn}.

The new dynamical system, \refe{fifodyn}, has the following properties.
(i) If $f_k^{rs}(0)=0$, then $f_k^{rs}(\tau)=0$ for any $\tau>0$; that is, if a route is unused initially, it is unused at any decision step.
(ii) The dynamical system is feasible in the following senses. First, if initial route flows are non-negative, they are always non-negative at any $\tau$, since the trajectories of $f_k^{rs}(\tau)$ cannot cross the boundary $f_k^{rs}=0$. Second, O-D flows are conserved; i.e., $\sum_k f_k^{rs}(\tau)=q_{rs}$, since $\sum_k \dot f_k^{rs}(\tau)=0$ for any $\tau$.

\begin{definition}[FIFO Equilibrium] Roots of equations $\bf J({\bf f})=\bf 0$ are {\em equilibria} of the dynamical system \refe{fifodyn}. Since O-D FIFO violation functions equal zero in these equilibria, an equilibrium of the dynamical system is also called a FIFO equilibrium.
\end{definition}

{\bf Remarks}. We have the following properties of a FIFO equilibrium.
\ben
\item ${\bf f}\in \mathcal{F}$ is an equilibrium of \refe{fifodyn} if and only if it satisfies one of the following conditions:
\ben
\item $J_k^{rs}=f_k^{rs}\sum_j (c_k^{rs}-c_j^{rs})f_j^{rs}=0$ ($\forall k,r,s$),
\item $J_k^{rs}=q_{rs}^2\xi_k^{rs}\sum_j (c_k^{rs}-c_j^{rs})\xi_j^{rs}=0$ ($\forall k,r,s$),
where $\xi_k^{rs}=f_k^{rs}/q_{rs}$ is the proportion of flow on route $k$ among flow on O-D pair $r-s$.
\item $f_k^{rs}(c_k^{rs}-v_{rs})=0$ ($\forall k,r,s$),
\item $c_k^{rs}=c_j^{rs}=v_{rs} \m{ when } f_k^{rs} f_j^{rs}\neq 0$ ($\forall k,j,r,s$),
\item $(c_k^{rs}-c_j^{rs})^n f_k^{rs} f_j^{rs}=0$ ($\forall k,r,s$), for any $n=1,2,\cdots$
\item $c_1^{rs}\leq \cdots \leq c_j^{rs} \leq c_{j+1}^{rs} =\cdots= c_{j+l}^{rs} \leq c_{j+l+1}^{rs} \leq \cdots$, where $f_k^{rs}>0, \forall\:\: k=j+1,\cdots, j+l$, and $f_k^{rs}=0, \forall\:\: k=1,\cdots, j, j+l+1,\cdots$.
\een
\item The dynamic system usually has multiple equilibria. For example, any vertex in the convex polygon defined by \refe{conservation} is an equilibrium. In equilibria, used routes share the same travel time, but unused routes may have longer or shorter travel times.
\item It is straightforward that, if ${\bf f}$ is a UE, then it is a FIFO equilibrium. However, some equilibria may not be UE. We call such FIFO equilibrium as {\em partial user equilibrium (PUE)}, since it can be considered as UE for a partial set of routes. We can see that a PUE satisfies the first part of Wardrop's definition of UE; i.e., in a PUE, ``the journey times of all routes actually used are equal". In real-world transportation networks, such PUE states could be realized if drivers prefer familiar routes over shorter ones connecting the same O-D pair.
\een

\subsection{Stability properties of UE and PUE}
In the following, we show that UE and PUE have different stability properties. First, we prove the instability of PUE.

\begin{theorem}\label{pue:unstable} A PUE is unstable for the dynamical system \refe{fifodyn}.
\end{theorem}
{\em Proof.} If ${\bf f}$ is a PUE, there exists an unused, shorter route $k$ for some O-D pair $r-s$. Then, for a used route $j$, $c_k^{rs}<c_j^{rs}$, $f_k^{rs}=0$, and $f_j^{rs}>0$. We introduce a perturbation to this equilibrium by shifting $\epsilon_0$ ($f_j^{rs}>\epsilon_0>0$) from route $j$ to route $k$. We denote the resultant state by $\tilde {{\bf f}}$ and the corresponding vector of route travel times by $\tilde {\bf c}$. Since route travel times are continuous in ${\bf f}$, we can have a sufficiently small $\epsilon_0$ such that route $k$ is still shorter than other used routes; i.e., $\tilde c_k^{rs}<\tilde c_j^{rs}$ and $\tilde c_k^{rs}<\tilde c_l^{rs}$ for any $f_l^{rs}>0$. Then at $\tilde f_k^{rs}=\epsilon_0$, the local changing direction in $\tilde f_k^{rs}$  is
\bqs
\dot {\tilde f_k^{rs}}&=&-\tilde f_k^{rs}\sum_l (\tilde c_k^{rs}-\tilde c_l^{rs}) \tilde f_l^{rs}>0.
\eqs
That is, the perturbed state $\tilde {\bf f}$ will drift away from the original equilibrium, and the original PUE is not stable. \eop

Then, we introduce a Lyapunov function for the dynamical system \refe{fifodyn} as follows.
\begin{theorem} \label{theorem:lya}
The function
\bqn
V({\bf f})=z({\bf f})-\min z({\bf f}) \label{lyapunov}
\eqn
is a Lyapunov function \citep{LaSalle1960lyapunov,Strogatz1994Nonlinear} of the dynamical system at UE. That is
\ben
\item $V({\bf f})>0$, $\forall\:\: {\bf f}\in \mathcal{F}$ and ${\bf f}$ is not a UE;
\item $V({\bf f})=0$ if and only if ${\bf f}\in\mathcal{E}$, where $\mathcal{E}$ is the set of UE;
\item -$grad \:V({\bf f})\cdot \J({\bf f})\leq 0$ if ${\bf f}$ is not a UE.
\een
This is equivalent to saying that the BMW objective function can be considered the potential energy of the dynamical system \refe{fifodyn}.
\end{theorem}
{\em Proof.} From the definitions of the Lyapunov function $V({\bf f})$ in \refe{lyapunov} and the BMW objective function $z({\bf f})$ in \refe{objective}, we can see that the first two statements are correct. We prove the third statement as follows.

We first compute the gradient of $V({\bf f})$ with respect to ${\bf f}$, whose $(rs,k)$th element is ($\forall\:\: k,r,s$)
\bqs
\pd {V({\bf f})}{f_k^{rs}}&=&\pd {z({\bf x})}{f_k^{rs}}=\sum_a t_a(x_a) \delta_{a,k}^{rs}=c_k^{rs}.
\eqs
Therefore,
\bqs
-grad \: V({\bf f})\cdot \J({\bf f})&=&\sum_{rs}\sum_k \pd {V}{f_k^{rs}} \dot f_k^{rs}=- \sum_{rs}\sum_k c_k^{rs} f_k^{rs}\sum_j (c_k^{rs}-c_j^{rs})f_j^{rs},
\eqs
which leads to
\bqn
-grad \: V({\bf f})\cdot \J({\bf f})&=&- \sum_{rs}\sum_k \sum_{j>k} (c_k^{rs}-c_j^{rs})^2 f_k^{rs} f_j^{rs}\leq 0,
\eqn
since all route flows are non-negative.
From the definition of equilibria, we can see that $-grad \: V({\bf f})\cdot \J({\bf f})=0$ if and only if ${\bf f}$ is an equilibrium, and $-grad \: V({\bf f})\cdot \J({\bf f})<0$ for non-equilibrium states. Therefore, $V({\bf f})$ is a Lyapunov function of the dynamical system \refe{fifodyn}. \eop

We have the following property for the set of equilibria of the dynamical system \refe{fifodyn}.
\begin{lemma}\label{lemma1}
The set of equilibria is closed. That is, given a sequence ${\bf f}_i\to {\bf f}^*$, if all ${\bf f}_i$ are equilibria, then ${\bf f}^*$ is also an equilibrium.
\end{lemma}
{\em Proof.}  Since ($\forall\:\: k,r,s$) \[c_k^{rs}({\bf f})=\sum_a \delta_{a,k}^{rs} t_a(\sum_{mn}\sum_l \delta_{a,l}^{mn} f_l^{mn})\] is continuous in ${\bf f}$, FIFO violation functions, $\J({\bf f})$, are also continuous in ${\bf f}$. Then, with ${\bf f}_i\to {\bf f}^*$, $\J({\bf f}_i)=0$ implies $\J({\bf f}^*)=0$. That is, ${\bf f}^*$ is also an equilibrium, and the set of equilibria is closed. However, note that the set of equilibria is generally not connected, since we may have multiple isolated regions of equilibria.  \eop

Further, we have the following properties for the set of UE.
\begin{lemma}\label{lemma2}
The set of UE, $\mathcal{E}$, is closed and connected.
\end{lemma}
{\em Proof.} Given a sequence of ${\bf f}_i\in \mathcal{E}$ and ${\bf f}_i\to {\bf f}^*$, then from continuity of the BMW objective function, we can see that $z({\bf f}_i)\to z({\bf f}^*)$. Therefore, $z({\bf f}^*)$ is also minimum, and ${\bf f}^*\in \mathcal{E}$. Thus $\mathcal{E}$ is closed. Since $z(\bf x)$ is a convex function and $\bf x$ is in a convex set, the set of UE, $\mathcal{E}$, is connected in $\bf x$. Further, since $x_a=\sum_{rs}\sum_k \delta_{a,k}^{rs} f_k^{rs}$ is a continuous mapping from $\mathcal{F}$ to $\mathcal{A}$, we have that $\mathcal{E}$ is also connected in ${\bf f}$.\eop

Finally, we show that UE is stable for symmetric link performance functions as follows.
\begin{theorem} \label{theorem:sta} A UE is locally, asymptotically stable for
the dynamical system \refe{fifodyn} in the sense that solutions of \refe{fifodyn} converge to a UE for a non-UE initial state that is close to the UE.
\end{theorem}
{\em Proof.} From the lemmas above, we can find a region around UE, where all non-UE states are not equilibria. Then $-grad \: V({\bf f}) \cdot \J({\bf f})<0$ for all the non-equilibrium states in this region. Therefore, according to the Lyapunov's stability theorem in \citep{smith1984stability}, asymptotic solutions of \refe{fifodyn} converge to UE solutions for non-UE initial states. That is, the set of UE is asymptotically stable. Note that we can have multiple UE and cannot use the traditional Lyapunov's stability test. \eop

From Theorems \ref{pue:unstable} and \ref{theorem:sta}, we can see that UE are the only stable equilibria of the new dynamical system for symmetric link performance functions. Note that, however, we can have unstable UE when link performance functions are non-monotone \citep{Netter1972critique}.

\subsection{A dynamical system model of static TAP with variable demand}
The dynamical system \refe{fifodyn} for fixed demand can be extended for static TAP with variable demand, whose mathematical programming formulation is given by \citep{Beckmann1956bmw,sheffi1985networks}
\bsq
\bqn
\min z({\bf x},{\bf q})&=&\sum_a \int_0^{x_a} t_a(\omega) d\omega-\sum_{rs}\int_0^{q_{rs}} u_{rs}(\omega)d\omega, \label{variabledemand}
\eqn
subject to \refe{conservation} and
\bqn
q_{rs}\geq0,&&\forall\:\: r,s,
\eqn
\esq
where $u_{rs}(q_{rs})$ is a decreasing travel time function associated with O-D pair $r-s$ and the inverse of the demand function. For symmetric link performance functions, the BMW objective function in \refe{variabledemand} is still convex, and there exists a unique user equilibrium in ${\bf x}$ and ${\bf q}$.

For this problem, we introduce the following dynamical system
\bsq\label{jsystem:variable}
\bqn
\dot f_k^{rs}&=&-q_{rs}f_k^{rs}(c_k^{rs}-u_{rs}), \quad \forall\:\:k,r,s\label{js1.1}\\
\dot q_{rs}&=&-q_{rs}(\sum_kf_k^{rs}c_k^{rs}-q_{rs}u_{rs}),\quad \forall\:\:r,s.\label{js1.2}
\eqn
\esq
Or, if introducing the average O-D travel time $v_{rs}$ as in \refe{averageOD}, we can have an equivalent dynamical system
\bsq\label{jsystem:variable2}
\bqn
\dot f_k^{rs}&=&-f_k^{rs}\sum_j(c_k^{rs}-c_j^{rs})f_j^{rs}-q_{rs}f_k^{rs}(v_{rs}-u_{rs}), \quad \forall\:\:k,r,s\label{js2.1}\\
\dot q_{rs}&=&-q_{rs}^2(v_{rs}-u_{rs}),\quad \forall\:\:r,s.\label{js2.2}
\eqn
\esq
Similar to \refe{fifodyn}, this dynamical system can have both UE and PUE as its equilibria. Further, we can have that: (i) PUE are unstable; (ii) the BMW objective function in \refe{variabledemand} can be considered as the potential energy function of the dynamical system \refe{jsystem:variable}; and (iii) only user equilibria are stable. These statements can be proved in the similar fashions as in the preceding subsection.

\section{Computation of static UE}
In the preceding section, we have derived a physically meaningful model \refe{fifodyn} for static TAP and analyzed its stability properties.
In this section, we introduce a numerical method for computing UE based on the dynamical system formulation. The purpose here is to show that the new formulation is feasible for finding UE numerically. Therefore, we do not concern with the computational efficiency of this method, which is an important factor for solving large-scale problems but not essential for analyzing median-sized prototype networks.

\subsection{Algorithms}
We first discuss algorithms for solving UE of static TAP with fixed demand.
Since an unused route in the initial state by vehicles of an O-D pair is always unused, we only consider those paths with non-zero initial flows and denote the number of initially used routes connecting O-D pair $r-s$ by $K_{rs}$ ($\forall\:\: r,s$).

By discretizing the dynamical system \refe{fifodyn} in the decision space as
\bqs
\frac{{\bf f}(\tau +\Delta \tau)-{\bf f}(\tau)}{\Delta \tau}=-\J({\bf f}(\tau)),
\eqs
we can update ${\bf f}$ by
\bqn
{\bf f}(\tau +\Delta \tau)={\bf f}(\tau)-\Delta \tau \J({\bf f}(\tau)),\label{fd}
\eqn
which is the standard Euler's method for solving ordinary differential equations \citep{Strogatz1994Nonlinear}.

Note that, at each step, we need to re-assign intermediate route flows to links to compute route travel times and obtain FIFO violation functions $\J({\bf f})$. Therefore, the total computational load of the finite difference method is proportional to the number of steps and that of computing FIFO violation functions. The index of convergence can be simply defined as the norm of FIFO violation functions as follows
\bqn
\| \J({\bf f})\|_2 &=&\sqrt{\sum_{rs}\sum_k (J_k^{rs})^2/\sum_{rs}K_{rs}}, \label{Jindex}
\eqn
which equals to zero if and only if ${\bf f}$ is an equilibrium.

For an arbitrary initial state ${\bf f}(0)$, the equilibrium of the dynamical system \refe{fifodyn} may not be UE, but a PUE, in which we can find unused shorter routes. For a PUE, we can perturb it by shifting proportions of flow from a used route to the shorter unused routes and obtain another equilibrium, since a PUE is unstable in these perturbation directions, as shown in the proof of Theorem \ref{pue:unstable}. We repeat this process until we reach a UE, where there are no shorter unused routes. The flow-chart of the algorithm is given in \reft{perturbation}. Note that the initial guess does not have to be an all-or-nothing assignment as in Frank-Wolfe's (FW) method \citep{sheffi1985networks}. Rather, with $K$-shortest routes as initial state, we could find UE with few or no perturbations.

For static TAP with variable demand, corresponding to the two formulations \refe{jsystem:variable} and \refe{jsystem:variable2} we can have two types of solution methods. In one, with an initial guess of $q_{rs}$ ($\forall r,s$) and $f_k^{rs}$ ($\forall k,r,s$), we directly solve \refe{js1.1} and \refe{js1.2}, by using the aforementioned perturbation-based method. In the other, with an initial guess of $q_{rs}$, we solve the corresponding UE $f_k^{rs}$ for a fixed demand problem \refe{js2.1} with $v_{rs}=u_{rs}$, and then update $q_{rs}$ and $f_k^{rs}$ with  the difference between $v_{rs}$ and $u_{rs}$ by \refe{js2.1} and \refe{js2.2} respectively.
\subsection{An example}
In this subsection, we study TAP on a simple network given by \citep[Fig. 5.1 on page 114 of][]{sheffi1985networks}, which is shown in \reff{network1}. For this network, FIFO violation functions can be written as
\bqs
J_1=x_1\left( \left(10(1+0.15(\frac{x_1}2)^4)-20(1+0.15(\frac{x_2}4)^4) \right)x_2+ \left( 10(1+0.15(\frac{x_1}2)^4)-25(1+0.15(\frac{x_3}3)^4) \right)x_3\right)\\
J_2=x_2\left( \left(20(1+0.15(\frac{x_2}4)^4)- 10(1+0.15(\frac{x_1}2)^4)\right)x_1+ \left(20(1+0.15(\frac{x_2}4)^4)-25(1+0.15(\frac{x_3}3)^4) \right)x_3\right)\\
J_3=x_3\left( \left(25(1+0.15(\frac{x_3}3)^4)- 10(1+0.15(\frac{x_1}2)^4)\right)x_1+ \left( 25(1+0.15(\frac{x_3}3)^4)-20(1+0.15(\frac{x_2}4)^4)\right)x_2\right)
\eqs
 and the dynamical system of the network is ($i=1,2,3$)
\bqs
\dot x_i=-J_i.
\eqs

First, we consider an initial state at $(x_1,x_2,x_3)=(3.39, 5.00, 1.61)$ and the corresponding route costs $(c_1,c_2,c_3)=(22.3, 27.3, 35.3)$. In the region of $x_1+x_2\leq 10$ ($x_3=10-x_1-x_2)$, the convergence is shown in \reff{fd:conv}, and the solution trajectories in \reff{fd:traj} with $\tau=0.02$ and $\Delta \tau=$0.0005, 0.001, and 0.002, respectively. From these figures, we can see that convergence rates and solution trajectories are  similar for different decision step sizes, as long as the finite difference form of \refe{fifodyn} is stable (i.e. $\Delta \tau$ is sufficiently small). In addition, convergence rates are constant when solutions get closer to equilibria, and there is no zigzagging effect.

With different initial conditions, we can find all seven equilibria for this network as shown in \reft{equilibria}.
We perturb the first six PUE by shifting $0.05$ flow from a used route to each unused shorter route. After perturbation, we use the finite difference method with $\Delta\tau=0.0005$ and $\tau=0.1$ to obtain solution trajectories of the dynamical system \refe{fifodyn}. As  shown in \reff{fig:pert}, all solution trajectories converge to UE after perturbation. This example confirms the theory that only UE are stable and demonstrates the feasibility of the perturbation-based method for finding UE.

\section{A dynamical system model of dynamic traffic assignment problem}
In this section, we consider dynamic TAP, in which all quantities are time-dependent. Since the concepts and definitions are all similar to those in static case, we do not specifically include ``dynamic" in them and keep the same terminologies as before. For example, UE is meant to be dynamic UE. We adopt the network notation system given in \reft{notation.3} and have the following remarks. First, flow-rates are derivatives of flows, and their units are number of vehicle per unit time and number of vehicles respectively. Second, all variables depend on the decision variable $\tau$, although not explicitly indicated. Third, the formulation is route-based, and we do not use link flows, link travel times, or indicator variables. In addition, we have the following basic relationships:
\bsq\label{conservation.3}
\bqn
p_{rs}(r,t)=\sum_k f_k^{rs}(r,t), \quad p_{rs}(s,t)=\sum_k f_k^{rs}(s,t),&& \forall\:\: r,s \\
f_k^{rs}(r,t)\geq 0, \quad f_k^{rs}(s,t)\geq 0,&& \forall\:\: k,r,s;
\eqn
where $p_{rs}$ and $f_k^{rs}$ are cumulative arrival curves at origins and destinations \citep{moskowitz1963capacity,newell1993sim},
and
\bqn
q_{rs}(r,t)=\sum_k g_k^{rs}(r,t), \quad q_{rs}(s,t)=\sum_k g_k^{rs}(s,t),&& \forall\:\: r,s \\
g_k^{rs}(r,t)\geq 0, \quad g_k^{rs}(s,t)\geq 0,&& \forall\:\: k,r,s.
\eqn
\esq

\subsection{Formulation}
For dynamic TAP, we propose the following dynamical system for adjusting route flow-rates
\bqn
-\der{g_k^{rs}(r,t)}{\tau}=J_k^{rs}(r,t), &&\forall\:\: k,r,s \label{dfds:con}
\eqn
where the time-dependent O-D FIFO violation function is defined by
\bsq\label{dynFV}
\bqn
J_k^{rs}(r,t)= q_{rs}(r,t)g_k^{rs}(r,t) (c_k^{rs}(r,t) - v_{rs}(r,t)),&&\forall\:\: k,r,s \label{cfvio:1}
\eqn
and the O-D average travel time for vehicles departing origin $r$ at time $t$, $v_{rs}(r,t)$, is
\bqn
v_{rs}(r,t)&=&\sum_j g_j^{rs}(r,t) c_j^{rs}(r,t)/q_{rs}(r,t).
\eqn
\esq
Note that \refe{dfds:con} is equivalent to $\der{g_k^{rs}(r,t)}{\tau}=-J_k^{rs}(r,t)$.
This dynamical system, continuous in time, is a direct extension of \refe{fifodyn} at any time $t$. Although there unlikely exist link performance functions when considering capacity constraints and link interactions \citep{daganzo1995tt}, the route travel time $c_k^{rs}(r,t)$ can be computed from a model of network vehicular traffic. For example, we can use the commodity-based kinematic wave model of network traffic developed in \citep{jin2003dissertation,jin2004network}, in which FIFO violation among vehicles on the same link converges to zero with diminishing cell sizes \citep{jin2005fifo}.

If O-D FIFO violation functions equal zero for all $r,s,k,t$, the corresponding state of flow-rates $g_k^{rs}(r,t)$ ($\forall\:\: k,r,s$) is called a dynamic {\em equilibrium}. We can have multiple equilibria, and trivial equilibria include those in which only one route is used for an O-D pair. Those equilibria are called dynamic {\em partial user equilibria} (PUE) if there are unused shorter routes,  and dynamic {\em user equilibria} (UE) otherwise. This definition of dynamic UE is consistent with Wardrop's. Further, from properties of static user equilibria, we have the following conjecture:
Dynamic PUE are unstable for \refe{dfds:con}; i.e., stable equilibria of \refe{dfds:con} are dynamic UE. In general, we expect to have multiple dynamic UE, and a dynamic UE can be unstable, due to non-monotonicity of link travel times with complicated interactions between traffic streams in a dynamic road network.

\subsection{Computation of dynamic UE}
In this subsection, we describe a computational method for finding equilibria and stable UE based on the new formulation \refe{dfds:con}. Here we assume that origin demands are given during an assignment time duration $[0, T_0]$ and zero outside the interval. We first discretize the assignment duration into $n$ time intervals with time instants $t_n=n\Delta t$, where $n=0,1,2,\cdots,N$ and $N=T_0/\Delta t$. We assume that during $[t_n,t_{n+1})$ the in-flow-rate of route $k$ connecting O-D pair $r-s$ is constant, $g_k^{rs}(r,t_n)$, and in-flows, $f_k^{rs}(r,t_n)$, are piece-wise linear functions.

In order to compute route travel times, we can use any traffic simulator with traffic simulation time duration of $[0,T]$ ($T>T_0$), during which all vehicles should be able to finish their trips. Also, in order to control numerical errors when computing travel times, we use smaller simulation time steps, $\Delta t/M$ ($M>1$), and simulation time instants of $t_i=i\frac{\Delta t}{M}$ ($i=0, 1, 2, \cdots, MN T/T_0$). Then, given an initial guess of $g_k^{rs}(r,t_n)$ ($\forall\:\: k,r,s$), boundary conditions at destinations, and initial traffic conditions, we can use a traffic simulator to obtain solutions of arrival curves at destinations $f_k^{rs}(s,t)$ at any simulation time instant.
To compute $c_k^{rs}(r,t_n)$, which is the average travel time of vehicles from $f_k^{rs}(r,t_n)$ to $f_k^{rs}(r,t_{n+1})$, we first find time instants $t_{m_n}$ and $t_{m_{n+1}}$ such that $f_k^{rs}(s,t_{m_n})$ and $f_k^{rs}(s,t_{m_{n+1}})$ are the closest to $f_k^{rs}(r,t_n)$ and $f_k^{rs}(r,t_{n+1})$ respectively. As illustrated in \reff{traveltime}, the total travel time of vehicles departing during $(n\Delta t, (n+1)\Delta t)$ is the area bounded by the two cumulative flow curves at the origin and destination $r-s$ and can be computed by
\bqn
c_k^{rs}(r,t_n)g_k^{rs}(r,t_n)\Delta t=\int_{f_k^{rs}(r,t_n)}^{f_k^{rs}(r,t_{n+1})} c_k^{rs}(r,t) d f_k^{rs}\approx (f_k^{rs}(r,t_{n+1}) - f_k^{rs}(r,t_n)) (\frac{m_{n+1}}{M}-n-\frac 12)\Delta t &&\nonumber \\-\sum_{i=m_n+1}^{m_{n+1}} (f_k^{rs}(s,t_i)- f_k^{rs}(s,t_{m_n}))\frac{\Delta t}{M}+ (f_k^{rs}(s,t_{m_{n+1}})- f_k^{rs}(s,t_{m_n})) \frac{\Delta t}{2M},&&\label{tt}
\eqn
from which we can compute the approximate O-D FIFO violation $J_k^{rs}(r,t_n)$ by using \refe{dynFV}. Then we solve \refe{dfds:con} by a finite difference approximation as follows
\bqn
g_k^{rs}(r,t_n)|_{\tau+\Delta \tau}= g_k^{rs}(r,t_n)|_{\tau} -J_k^{rs}(r,t_n) \Delta \tau,&&\forall\:\: k,r,s,
\eqn
where $n=0,1,\cdots, N-1$. From these solutions, we can construct a new time series of in-flow-rates, $g_k^{rs}(r,t_n)$, at $\tau+\Delta \tau$. Here we use the following the index to measure the convergence of solutions
\bqn
\|\J\|_2&=&\sqrt{\frac{\sum_{rs}\sum_{k} \sum_{n=1}^N (J_k^{rs}(r,t_n))^2}{N\sum_{rs} K_{rs}}}, \label{Jindex.3}
\eqn
where $K_{rs}$ is the number of initially used routes connecting O-D pair $r-s$. If the finite difference equations are stable and convergent for a sufficiently small $\Delta \tau$, solutions will converge to equilibria for a sufficiently large $\tau$. The flow-chart of this algorithm for finding an equilibrium is given in \reft{alg:equilibria}. Further, if there are shorter unused routes in a PUE, we shift some flows to such routes and find another equilibrium, until we reach a stable UE.

\subsection{An example}
We study dynamic TAP in the network shown in \reff{networkdta}, where two routes have the same triangular fundamental diagram \citep{munjal1971multilane,newell1993sim}. The network is empty initially, and there is no capacity constraint at the destination. Here we consider a constant arrival flow-rate at the origin, $q_0$, during time interval $[0,T_0]$. That is, $q_{rs}(r,t)=q_0$ for $t\in[0,T_0]$ and 0 otherwise.

We use the Lighthill-Whitham-Richards \citep{lighthill1955lwr,richards1956lwr} traffic model to analyze traffic dynamics on each link. For example, for link 1, however many vehicles are waiting to enter, the maximum flow-rate is $q_c=1$. Since there is no bottleneck on the link or the destination, we know that traffic on link 1 is free flow, and the travel time is 1. Although the link travel time is constant, the waiting time at the origin is time-dependent. If the arrival flow at the destination for route 1 is $f_1(s,t)$, then the departure flow at the origin is $f_1(s,1+t)$. Denoting $t_i=i\Delta t/M$, we then have the following equation
\bsq
\bqn
f_1(s,1+t_{i+1})&=&f_1(s,1+t_i)+\frac {\Delta t}{M}\min\{q_c, \frac{f_1(r,t_i)-f_1(s,1+t_i)}{\Delta t}M + g_1(r,t_i) \}, \label{model:link1}
\eqn
where $f_1(s,1)=0$, $f_1(r,0)=0$, and $\frac{f_1(r,t_i)-f_1(s,1+t_i)}{\Delta t}M + g_1(r,t_i)$ is the maximum flow-rate that can be sent from the origin if there is no capacity constraint on this link. This method of computing boundary fluxes is based on the supply-demand method developed in \citep{daganzo1995ctm,lebacque1996godunov,jin2003inhlwr}. Thus, \refe{model:link1} gives the dynamic model for finding arrival flows at the destination on route 1. Similarly, we can have the following model for route 2:
\bqn
f_2(s,2+t_{i+1})&=&f_2(s,2+t_i)+\frac {\Delta t}{M}\max\{q_c, \frac{f_2(r,t_i)-f_2(s,2+t_i)}{\Delta t}M + g_2(r,t_i) \},
\eqn
where $f_2(s,2)=0$ and $f_2(r,0)=0$.
\esq

We set $T_0=1$, $T=8$, $N=20$, $M=10$, $q_0=5$, $\tau=160$, and $\Delta \tau=0.05$, and consider deterministic initial time series of $g_1(r,t_n)=cq_0$ and $g_2(r,t_n)=(1-c)q_0$ ($n=1,\cdots,N$). When $c=0.5$; i.e., both routes have the same flows initially, solutions of cumulative arrival curves of route flows and total flows at the origin and the destination are shown in \reff{dtasolution}, from which we can clearly see that these solutions satisfy the O-D FIFO principle. At different time instants, the travel times are shown in \reff{dtatt}, from which confirms that the solution is an equilibrium given by
\bsq\label{analyticalsol}
\bqn
f_1(r,t)&=&\cas{{ll}5t,& 0\leq t\leq 0.25;\\ 1.25 +2.5(t-0.25),& 0.25<t\leq 1;\\3.125, & t>1,}\\
f_2(r,t)&=&\cas{{ll}0,& 0\leq t\leq 0.25;\\ 2.5(t-0.25),& 0.25<t\leq 1;\\1.875, & t>1.}
\eqn
\esq
Since there are no unused shorter route in the solution, the equilibrium in \refe{analyticalsol} is a UE.

With the aforementioned deterministic initial time series with $c=0.5$, 0.95, and 0.05, and another random initial time series, we find that all solutions converge to \refe{analyticalsol} and the index of convergence is shown in \reff{dtaconvergence}, from which we can see that solutions converge slowly at the beginning stage, then exponentially when it is close to \refe{analyticalsol}. Note that $g_1(r,t_n)=q_0$ ($\forall\:\: n=1,\cdots, N$) and $g_2(r,t_n)=q_0$ ($\forall\:\: n=1,\cdots, N$) are two trivial equilibria. Since initial time series with $c=0.95$ can be considered as a perturbation around $g_1(r,t)=q_0$ and those with $c=0.05$ a perturbation around $g_2(r,t)=q_0$, we can conclude that the UE in \refe{analyticalsol} is stable. Convergence with random initial time series further confirms the stability of the UE.

\section{Discussions}
In this study, based on a definition of O-D FIFO violation among different routes connecting the same O-D pair, we derived and analyzed a new dynamical system formulation of the traffic assignment problem for both static and dynamic road networks. For static, symmetric link performance functions, we also showed that the dynamical system is stable at UE, and the well-known BMW objective function can be considered the potential energy of the new dynamical system. Through this study, we can see that the new model is physically meaningful, theoretically rigorous, and numerically feasible.

The new dynamical system model is simpler in form than existing aggregate, deterministic dynamical systems proposed in \citep{smith1984stability,Nagurney1996ds,Nagurney1997equilibria}. Also different from existing dynamical systems, its equilibria can be UE or PUE. Compared with mathematical programming formulations, this formulation is applicable for non-monotone traffic assignment problem, and we can apply this dynamical system model to analyze the stability of multiple UE. In addition, this new approach can be extended for other traffic assignment problems. For example, since system-optimal assignments are equivalent to user-equilibrium assignments with modified link performance functions, we can also have a similar dynamical system formulation of system-optimal assignments. In addition to being applied to theoretical formulation and analysis of various traffic assignment problems, the new model can also be used to numerically study small- to medium-sized networks, which can be solved in a reasonable amount of time with modern computers.

In the future, we will be interested in further theoretical and numerical investigations of the new dynamical system. Theoretically, the existence of stable equilibria of the dynamical system of static TAP with monotone link performance functions is subject to further investigations; it is also possible to incorporate network control policies and drivers' familiarity and preference of certain routes into the dynamical system; it will also be interesting to study the existence, uniqueness, and stability properties of dynamic UE in a general network with this formulation.
Numerically, although the Euler-based perturbation method proposed in this study is shown to be feasible with two simple examples, the method only converges linearly and is not sufficiently efficient for large-scale problems. In the future, we will be interested in developing more efficient methods based on this formulation. For examples, we could develop parallel computational methods, or hybrid methods integrating the Frank-Wolfe method and the perturbation-based method; we can also use the fourth-order Runge-Kutta method, which is more efficient than the Euler's method.

\section*{Acknowledgement}
I greatly appreciate Hyunmyung Kim and R. Jayakrishnan for some intriguing discussions on the traffic assignment problem. I would like to appreciate Prof. Will Recker for providing me a friendly research environment. I would also like to thank Fred Mannering, Fan Yang, Ding Zhang, H Michael Zhang, and two anonymous referees for their useful comments and suggestions. The views and results contained herein are the author's alone.

\clearpage

\btb
\bc
\begin{tabular}{ll}\\\hline
$\mathcal{N}$&node (index) set\\
$\mathcal{A}$&arc (index) set\\
$\mathcal{R}$&set of origin nodes; $\mathcal{R}\subset\mathcal{N}$\\
$\mathcal{S}$&set of destination nodes; $\mathcal{S}\subset\mathcal{N}$\\
$\mathcal{K}_{rs}$&set of routes connecting O-D pair $r-s$; $r\in \mathcal{R}$, $s\in \mathcal{S}$\\
$\tau$ & decision variable, independent of time\\
$\Delta \tau$ & decision step\\
$x_a$ & flow on arc $a$; ${\bf x}=(\cdots, x_a,\cdots)$\\
$t_a$ & travel time on arc $a$;${\bf t}=(\cdots, t_a,\cdots)$\\
$q_{rs}$ & traffic demand between origin-pair $rs$; $q_{rs}\equiv f_0^{rs}$\\
$v_{rs}$& average travel time for O-D pair $r-s$\\
$u_{rs}$ & travel time function for variable demand for O-D pair $r-s$\\
$f_k^{rs}$ & flow on route $k$ connecting O-D pair $r-s$; ${\bf f}^{rs}=(f_0^{rs},\cdots, f_k^{rs}, \cdots)$; ${\bf f}=(\cdots, {\bf f}^{rs}, \cdots)$\\
$c_k^{rs}$ & travel time on route $k$ connecting O-D pair $r-s$; ${\bf c}^{rs}=(\cdots, c_k^{rs}, \cdots)$; $\bf c=(\cdots, \bf c^{rs}, \cdots)$\\
$\delta_{a,k}^{rs}$ & indicator variable: $\delta_{a,k}^{rs}=\cas{{ll} 1 & \m{if link } a \m{ is on route } k \m{ between O-D pair }r-s\\0&\m{otherwise;} } $ \\
&$\pmb \Delta^{rs}=(\cdots, \delta_{a,k}^{rs},\cdots)$; $\pmb \Delta=(\cdots, \pmb \Delta^{rs}, \cdots)$\\
$\mathcal{F}^{rs}$ & the set of ${\bf f}^{rs}$ satisfying \refe{conservation}\\
$\mathcal{F}$ & the set of ${\bf f}$; $\mathcal{F}=\prod_{rs} \mathcal{F}^{rs}$\\
$\mathcal{G}^{rs}$ & the set of initially non-empty ${\bf f}^{rs}$ satisfying \refe{conservation} and ${\bf f}^{rs}(\tau)>{\bf 0}$ for $\tau=0$\\
$\mathcal{G}$ & the set of ${\bf f}$; $\mathcal{G}=\prod_{rs} \mathcal{G}^{rs}$\\
$K_{rs}$ & the number of initially non-empty routes connecting O-D pair $r-s$\\
$J_k^{rs}$ & FIFO violation for flow on route $k$ connecting O-D pair $r-s$; \\ & ${\bf J}^{rs}=(\cdots, J_k^{rs}, \cdots)$; ${\bf J}=(\cdots, {\bf J}^{rs}, \cdots)$\\
$\|{\bf J}\|_2$ & 2-norm of ${\bf J}({\bf f})$, defined in \refe{Jindex}
\\\hline
\end{tabular}\caption{Network notation system for static TAP}\label{notation.2}
\ec
\etb

\bfg
\bc
\includegraphics{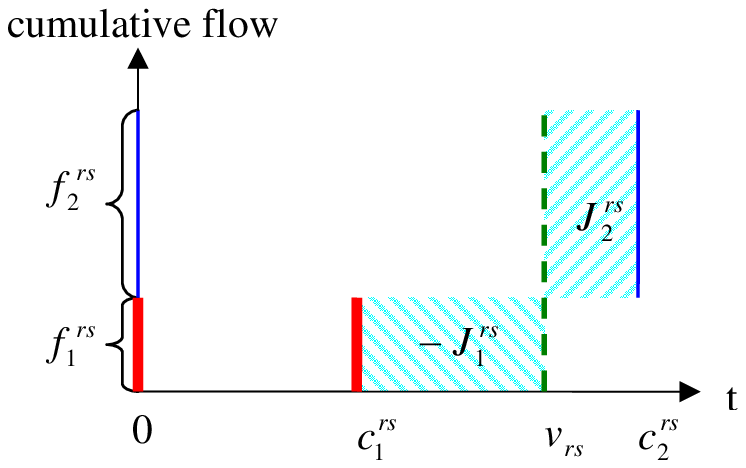}
\caption{An example of O-D First-In-First-Out violation}\label{fifostatic}
\ec
\efg

 \btb\bc
\begin{tabular}{|l|}\hline
Initial guess of ${\bf f}$, e.g. $K$-shortest routes for each O-D pair\\
for $n=1,2,3\cdots$\\
$\qquad$ Find an equilibrium with Euler's method\\
$\qquad$ If no shorter unused routes, the equilibrium is UE;\\
$\qquad$ Otherwise, perturb the equilibrium in the directions of shorter unused routes \\\hline
\end{tabular}\caption{A perturbation-based method for finding stable UE} \label{perturbation}
\ec\etb

\bfg\bc
\includegraphics{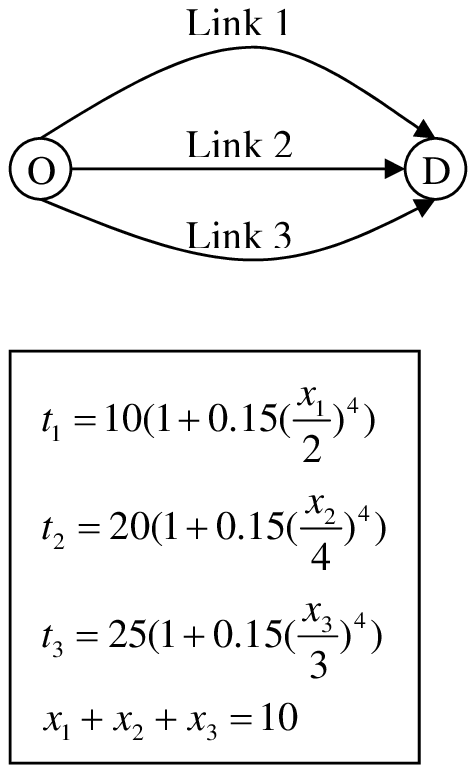}\caption{A static road network}\label{network1}
\ec\efg

\bfg\bc
\includegraphics[height=8cm]{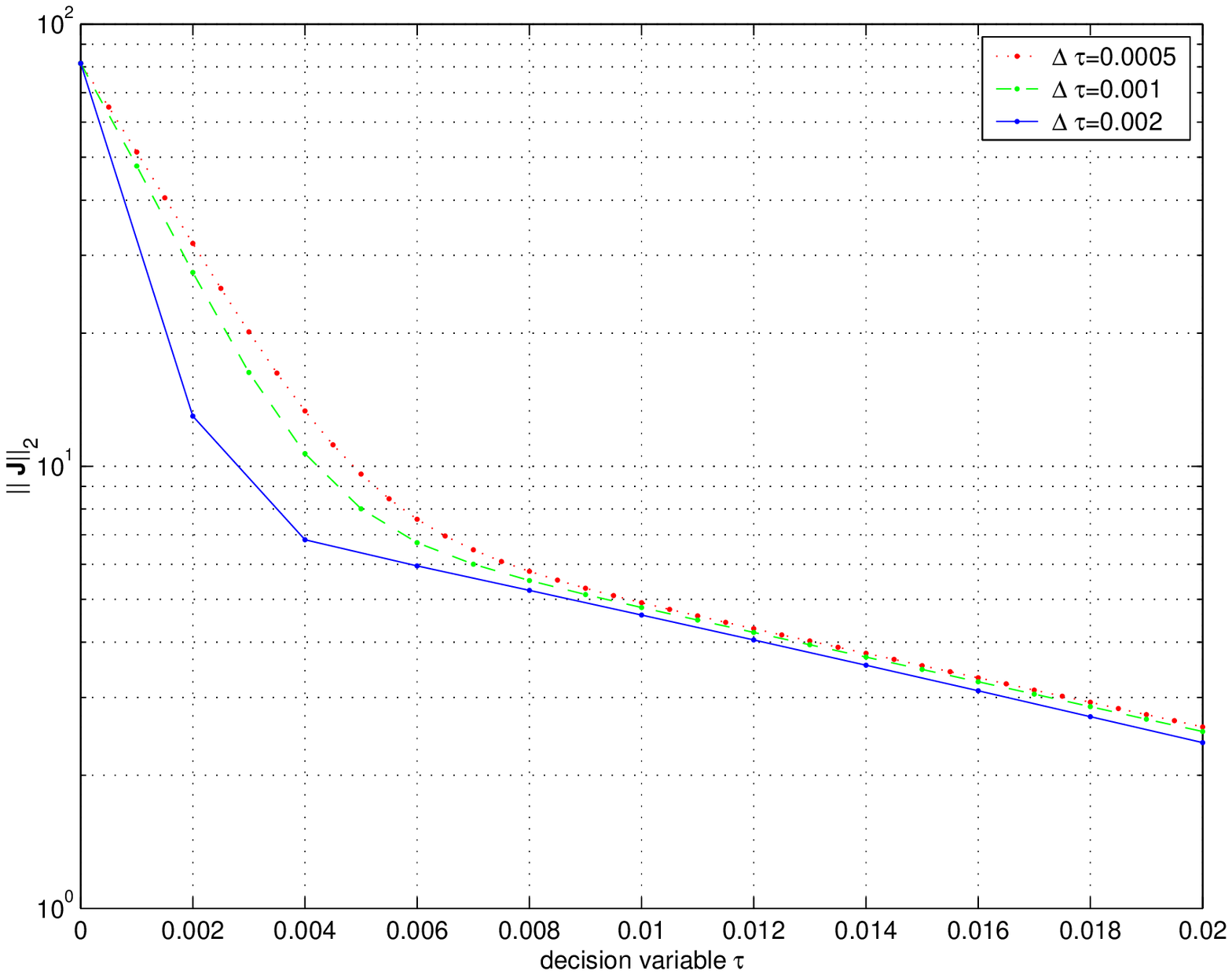}\caption{Convergence of Euler's method} \label{fd:conv}
\ec\efg

\bfg\bc
\includegraphics[height=8cm]{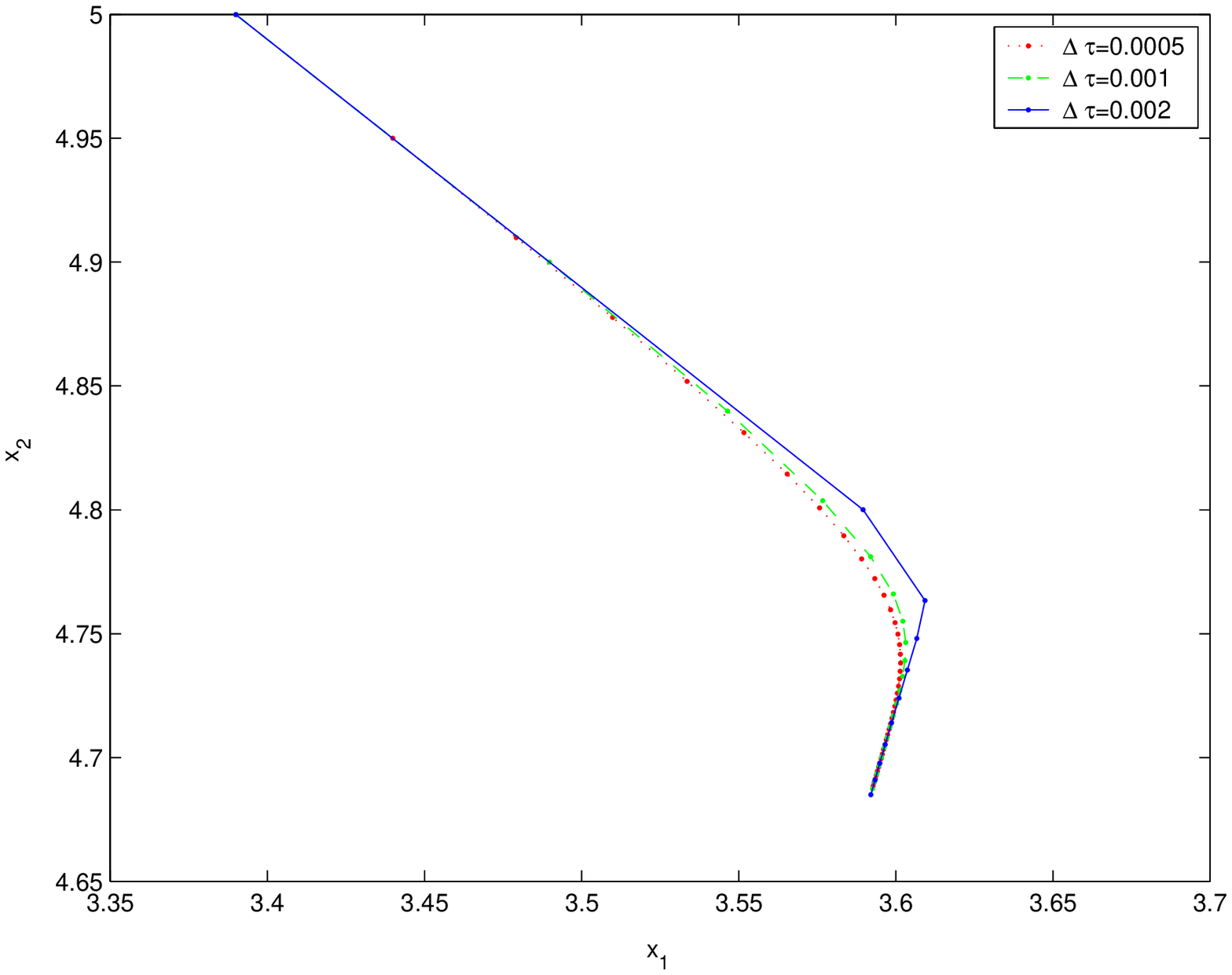}\caption{Solution trajectory of Euler's method with $\tau=0.02$} \label{fd:traj}
\ec\efg

\btb
\bc
\begin{tabular}{|c|ccc||ccc|}\hline
equilibria&\multicolumn{3}{|c||}{link flow}&\multicolumn{3}{|c|}{link cost}\\\hline
1&10& 0& 0 & 947.5000&   20  & 25\\
2&0&10&0&10 & 137.1875&25\\
3&0&0&10&10&20 & 487.9630\\
4&4.0346 &5.9654&0&34.8405&34.8405&25\\
5&4.7864&0 &5.2136&59.2053&20&59.2053\\
6& 0 &6.0762 &3.9238&10&35.9740&35.9740\\
7&3.5833 &4.6451 &1.7716&25.4560&25.4560&25.4560\\
\hline
\end{tabular}\caption{All equilibria for network in \reff{network1}}\label{equilibria}
\ec\etb

\bfg\bc
\includegraphics[height=8cm]{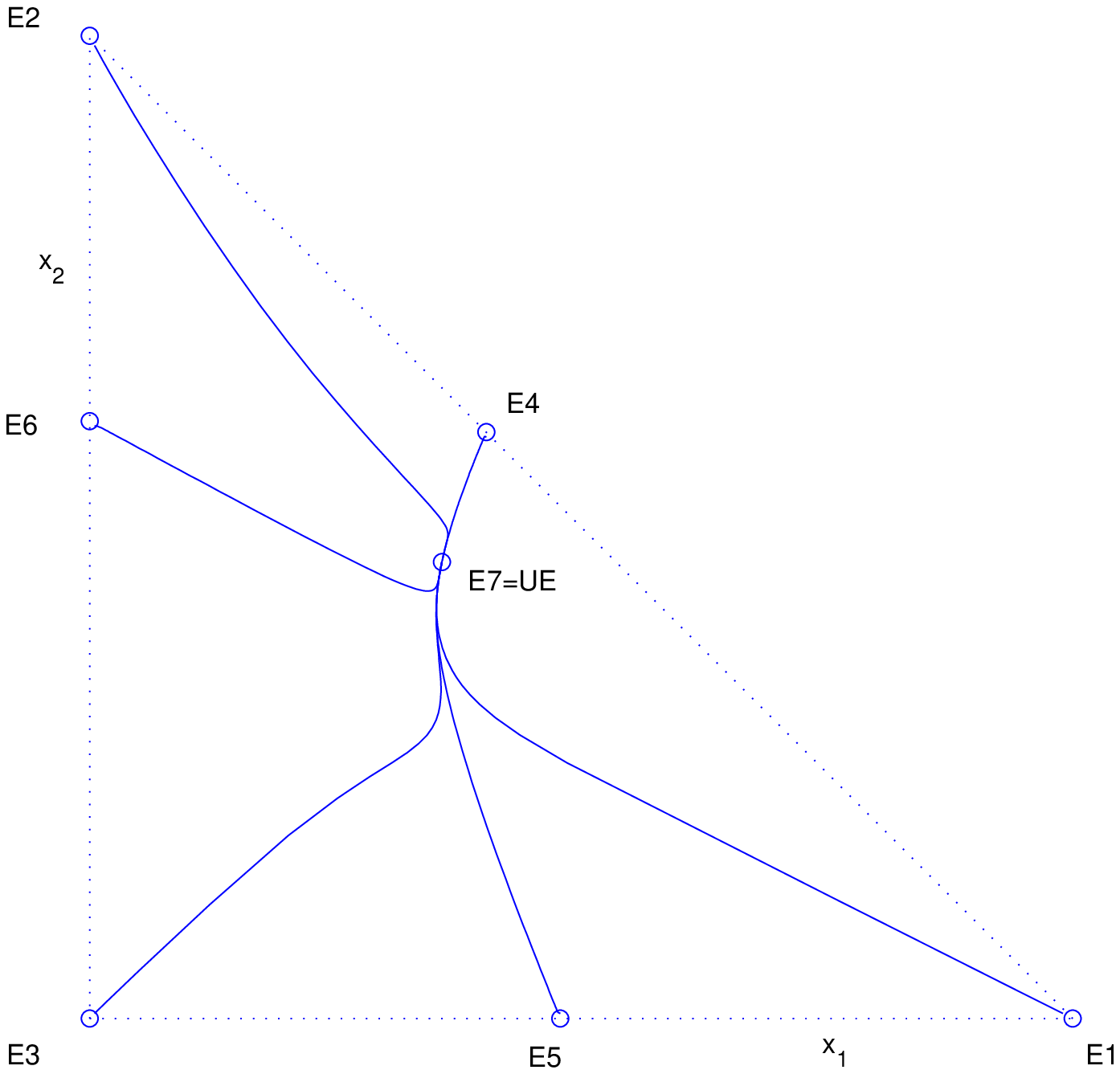}\caption{Stability of UE and instability of PUE} \label{fig:pert}
\ec\efg

\btb
\bc
\begin{tabular}{ll}\\\hline
$\mathcal{R}$&set of origin nodes; $\mathcal{R}\subset\mathcal{N}$\\
$\mathcal{S}$&set of destination nodes; $\mathcal{S}\subset\mathcal{N}$\\
$\mathcal{K}_{rs}$&set of routes connecting O-D pair $r-s$; $r\in \mathcal{R}$, $s\in \mathcal{S}$\\
$\tau$ & independent decision variable\\
$\Delta \tau$ & decision step\\
$t$ & independent time variable\\
$[0, T_0]$ & assignment interval\\
$[0, T]$ & traffic simulation interval\\
$N$ & the number of assignment time steps\\
$M$ & the number of traffic simulation steps during each assignment time step\\
$\Delta t$ & assignment time step, $\Delta t=T_0/N$\\
$\Delta t/M$ & traffic simulation time step\\
$t_n$ & assignment time instants, $t_n=n\Delta t$ for $n=0, \cdots N$\\
$t_i$ & traffic simulation time instants, $t_i=i \Delta t/M$ for $i=0,\cdots, MN T/T_0$\\
$f_k^{rs}(r,t)$ & arrival flow at origin $r$ at $t$ taking route $k$ connecting O-D pair $r-s$; \\
&${{\bf f}}^{rs}(r,t)=(\cdots, f_k^{rs}(r,t), \cdots)$; ${{\bf f}}(r,t)=(\cdots, {{\bf f}}^{rs}(r,t), \cdots)$\\
$g_k^{rs}(r,t)$ & arrival flow-rate at origin $r$ at $t$ taking route $k$ connecting O-D pair $r-s$;\\
&${\bf g}^{rs}(r,t)=(\cdots, g_k^{rs}(r,t), \cdots)$; ${\bf g}(r,t)=(\cdots, {\bf g}^{rs}(r,t), \cdots)$\\
$f_k^{rs}(s,t)$ & arrival flow at destination $s$ at $t$ taking route $k$ connecting O-D pair $r-s$; \\
$g_k^{rs}(s,t)$ & arrival flow-rate at destination $s$ at $t$ taking route $k$ connecting O-D pair $r-s$; \\
$p_{rs}(r,t)$ & arrival flow at origin $r$ between origin-pair $r-s$; \\
$q_{rs}(r,t)$ & arrival flow-rate at origin $r$ between origin-pair $r-s$; \\
$p_{rs}(s,t)$ & arrival flow at destination $s$ between origin-pair $r-s$; \\
$q_{rs}(s,t)$ & arrival flow-rate at destination $s$ between origin-pair $r-s$; \\
$c_k^{rs}(r,t)$ & travel time on route $k$ connecting O-D pair $r-s$ for a vehicle arriving at origin $r$ at $t$; \\
$v_{rs}(r,t)$ & travel time between O-D pair $r-s$ for a vehicle arriving at origin $r$ at $t$; \\
$\mathcal{F}^{rs}(t)$ & the set of ${{\bf f}}^{rs}(t)$ satisfying \refe{conservation.3}\\
$\mathcal{F}(t)$ & the set of ${{\bf f}}(t)$; $\mathcal{F}(t)=\prod_{rs} \mathcal{F}^{rs}(t)$\\
$K_{rs}$ & the number of initially non-empty routes connection O-D pair $r-s$\\
$J_k^{rs}(r,t)$ & FIFO violation of flow from destination $r$ at $t$ on route $k$ connecting O-D pair $r-s$; \\ & ${\bf J}^{rs}(r,t)=(\cdots, J_k^{rs}(r,t), \cdots)$; ${\bf J}(r,t)=(\cdots, {\bf v}^{rs}(r,t), \cdots)$; $\J=(\cdots, \J(r,t_n),\cdots )$\\
$\|\J\|_2$ & 2-norm of $\J({\bf f})$, defined in \refe{Jindex.3}
\\\hline
\end{tabular}\caption{Network notation system for dynamic TAP}\label{notation.3}
\ec
\etb

\bfg
\bc\includegraphics[height=8cm]{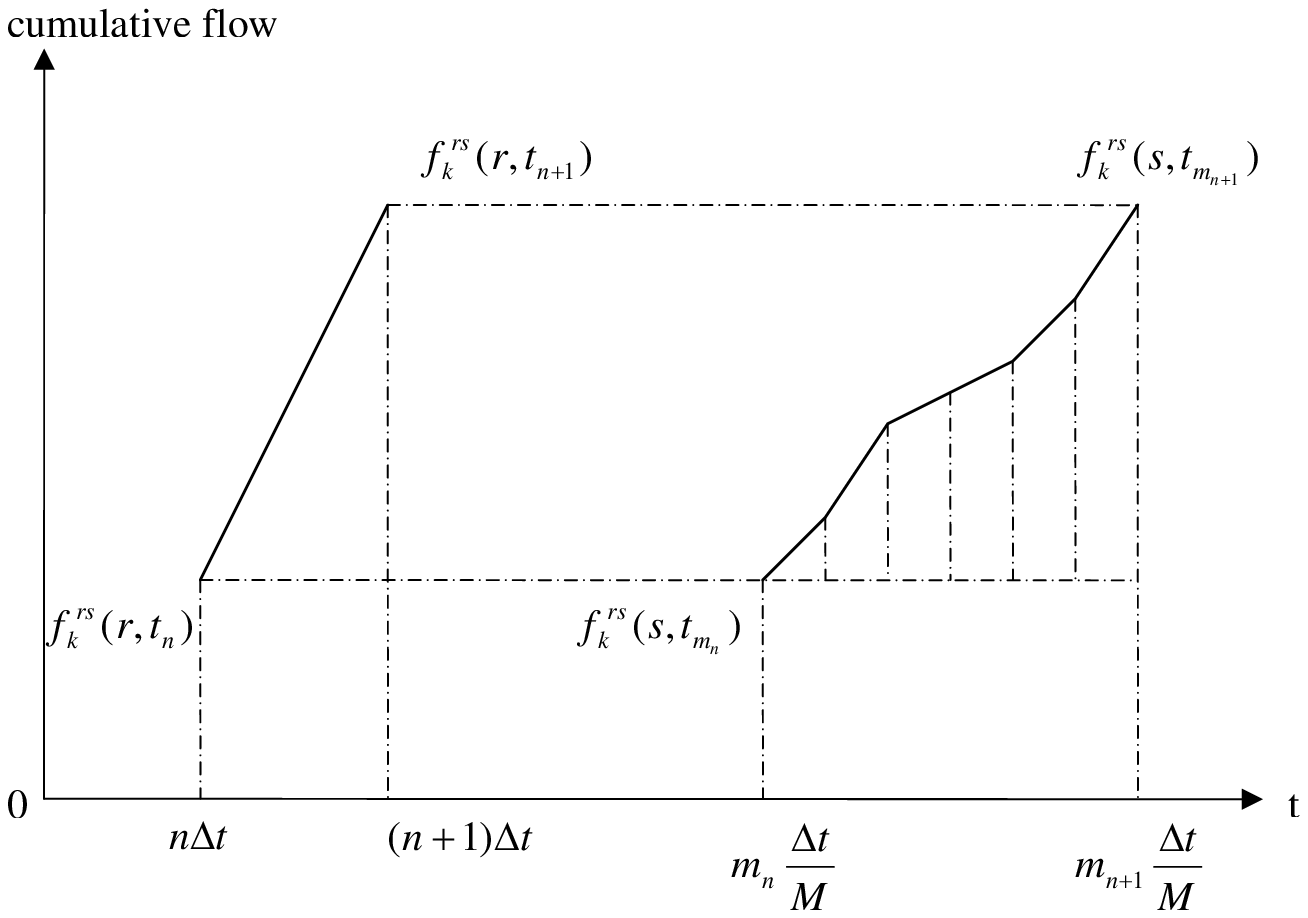}\caption{Computation of dynamic route travel time from simulated accumulative curves}\label{traveltime}\ec
\efg

\btb\bc
\begin{tabular}{|l|}\hline
Initialize $T_0$, $T$, $\tau$, $N$, $M$, $\Delta t$, $\Delta \tau$\\
Initial guesses of $g_k^{rs}(r,t_n)$ $\forall\:\: k,r,s,n$\\
for $\tau_l=\Delta \tau,2\Delta \tau,3\Delta \tau \cdots$\\
$\qquad$ Use traffic flow model to obtain $f_k^{rs}(s,t_i)$ $\forall\:\: k,r,s$ and $i=1,\cdots, MN T/T_0$\\
$\qquad$ Compute travel time $c_k^{rs}(r, t_n)$ and $v_{rs}(r, t_n)$\\
$\qquad$ Compute FIFO violation functions $J_k^{rs}(r,t_n)$ $\forall\:\: k,r,s$\\
$\qquad$ Use Euler's method to solve $g_k^{rs}(r,t_n)$ and find $f_k^{rs}(r,t_n)$ ($\forall\:\: n=1,\cdots,N$)\\
 \\\hline
\end{tabular}\caption{An algorithm for finding a dynamic equilibrium} \label{alg:equilibria}
\ec\etb

\bfg
\bc\includegraphics[height=8cm]{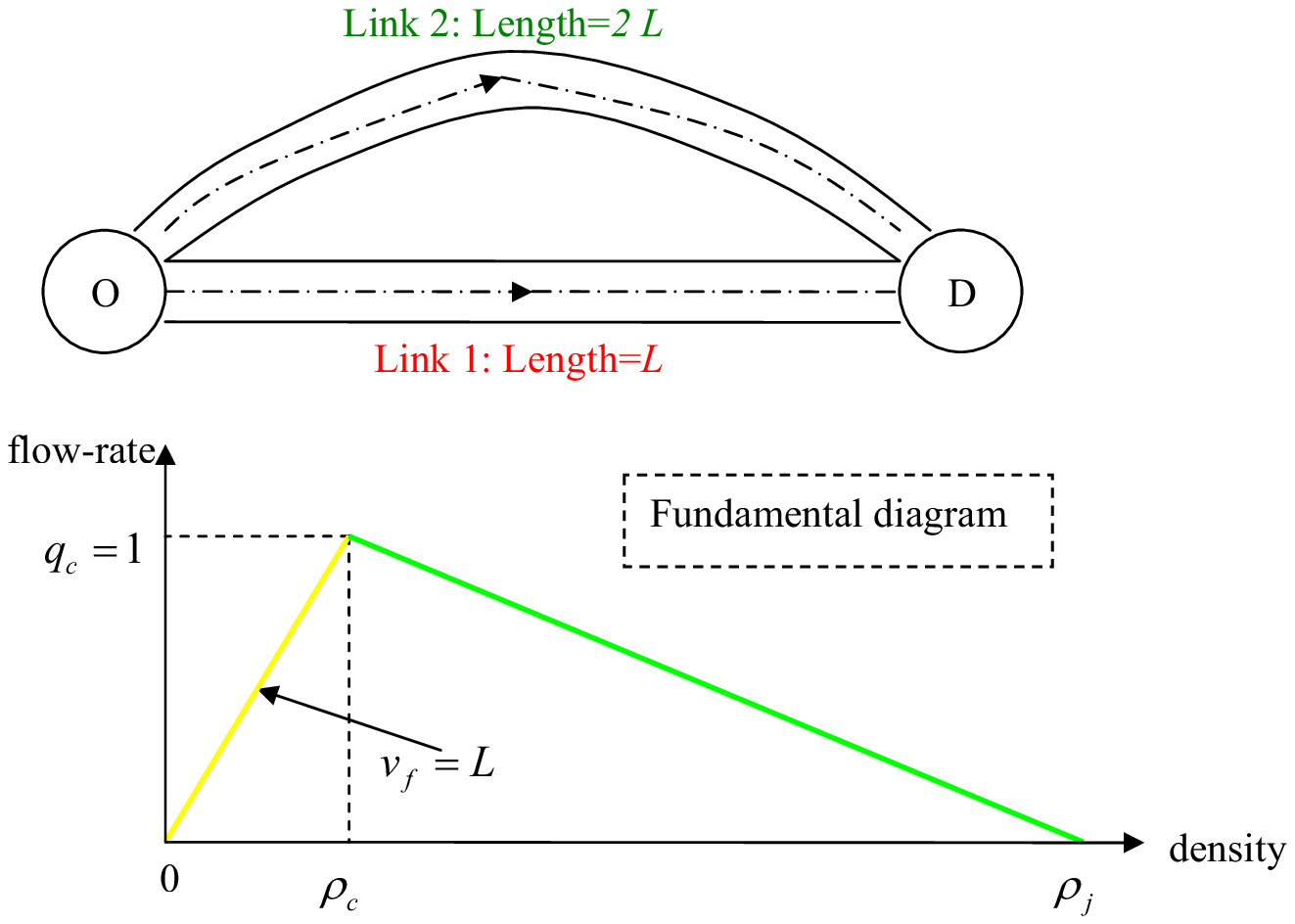}\caption{A dynamic road network}\label{networkdta}\ec
\efg

\bfg
\bc\includegraphics[height=8cm]{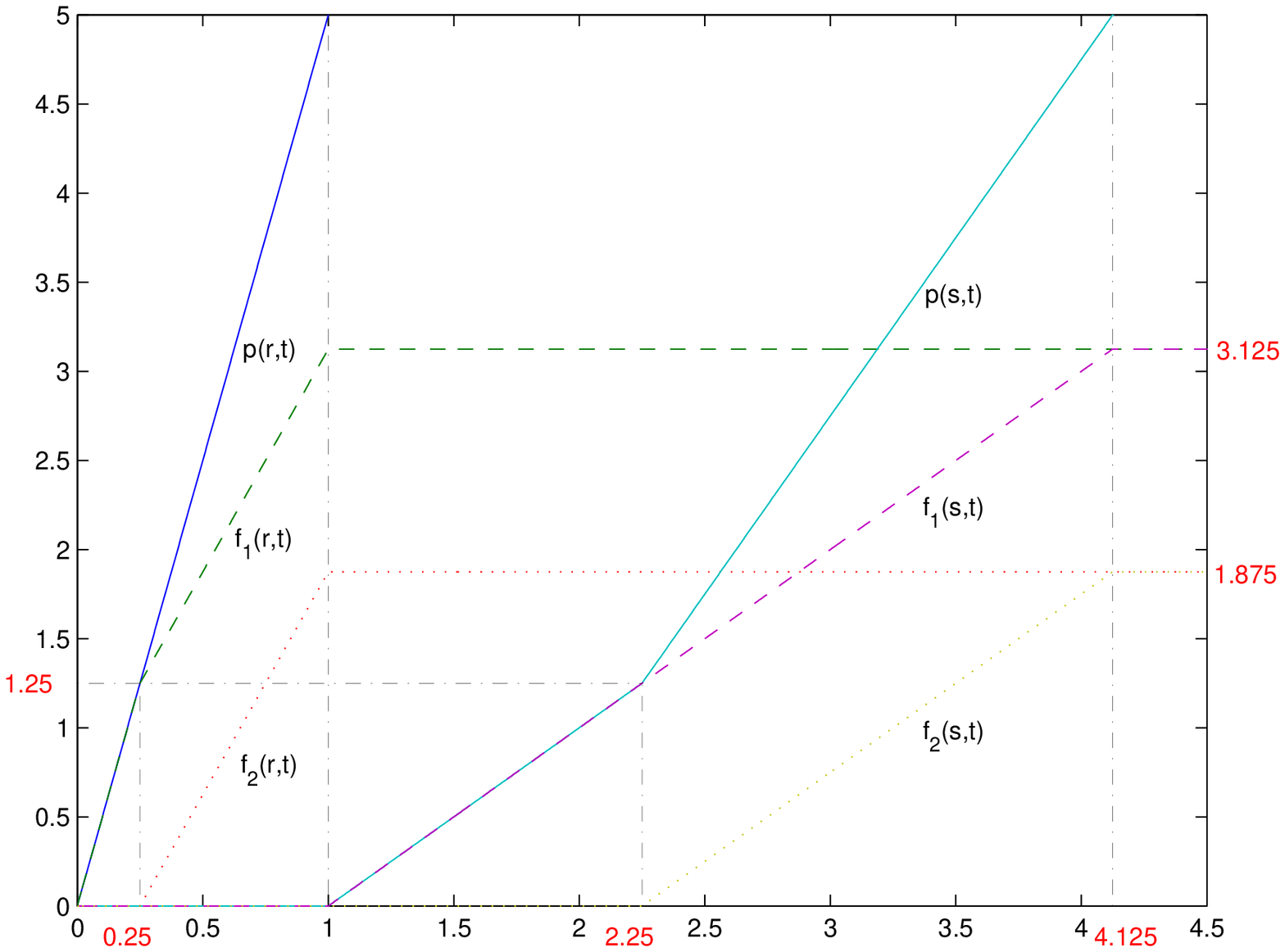}\caption{A dynamic UE for the network in \reff{networkdta}}\label{dtasolution}\ec
\efg

\bfg
\bc\includegraphics[height=8cm]{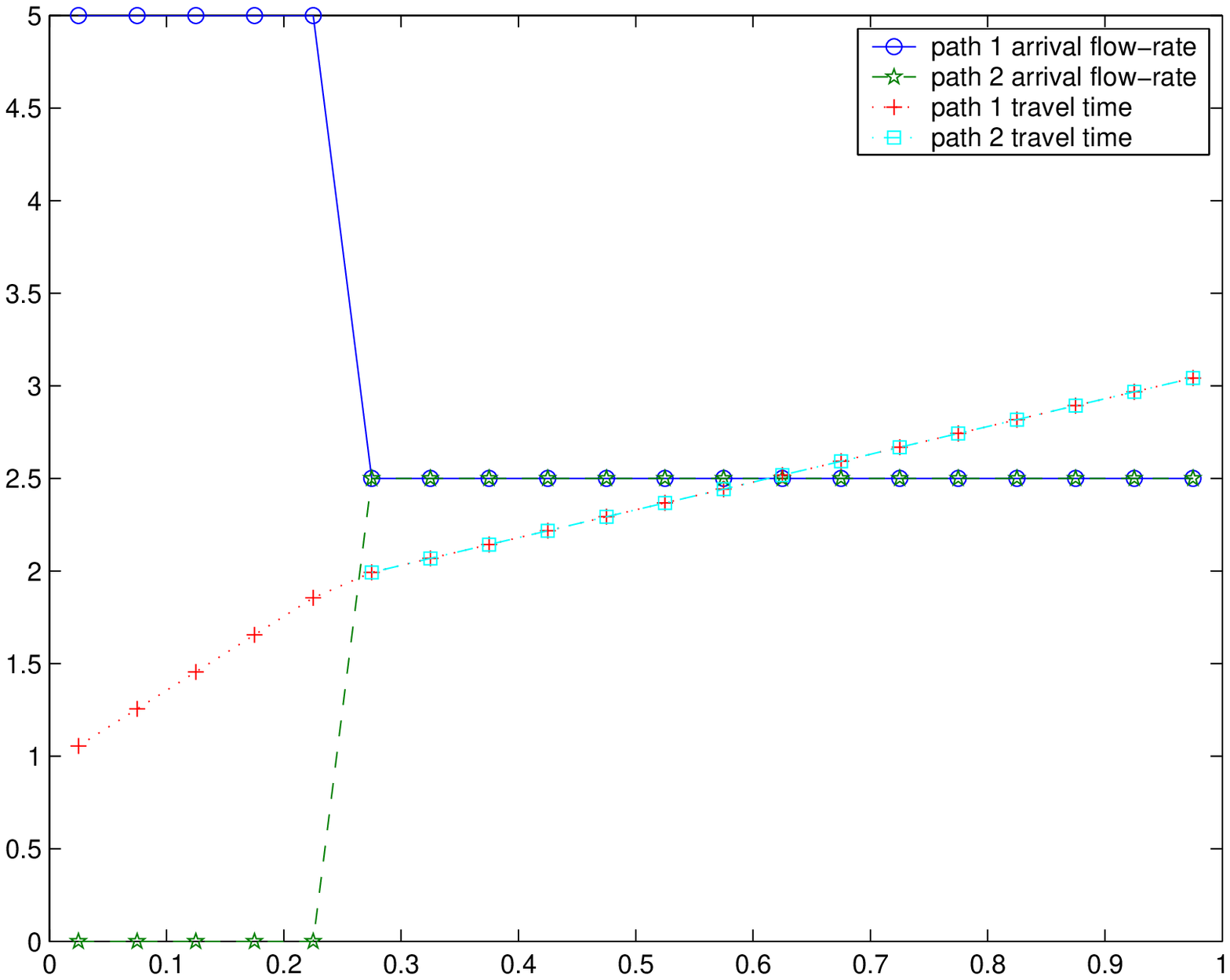}\caption{UE arrival flow-rates at the origin and travel times for the network in \reff{networkdta}}\label{dtatt}\ec
\efg

\bfg
\bc\includegraphics[height=8cm]{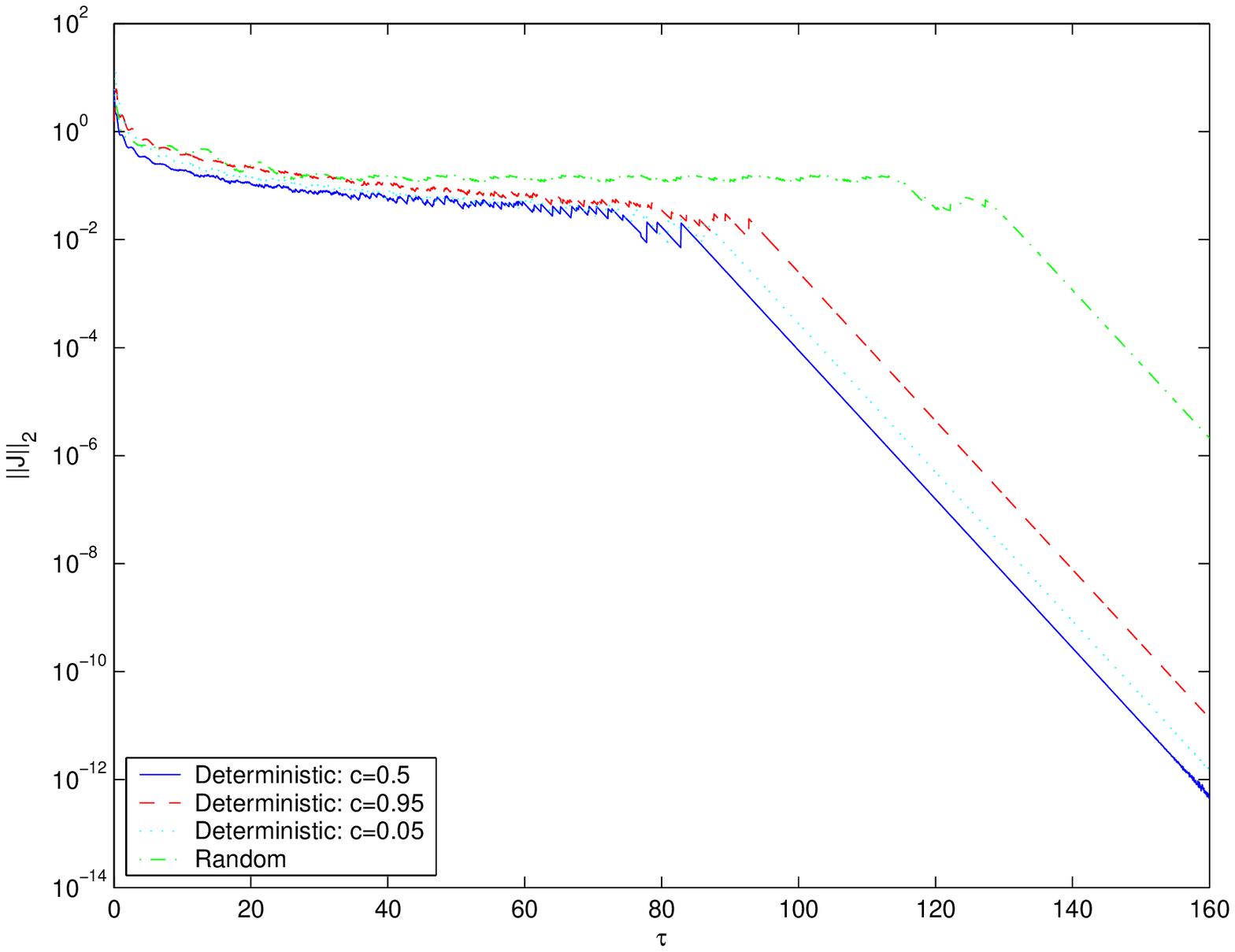}\caption{Convergence of solutions to a dynamic UE with different initial conditions with $\Delta \tau=0.05$ for the network in \reff{networkdta}}\label{dtaconvergence}\ec
\efg

\end{document}